\documentclass[11pt]{article}
\usepackage{amssymb}
\usepackage{amsthm}
\usepackage{amsmath}
\usepackage{fullpage}
\usepackage{epsfig}
\usepackage{hyperref}

\newcommand\remove[1]{}

\newcommand{\rnote}[1]{}

\newcommand{\f}{\varphi}
\newcommand{\Lip}{\mathrm{Lip}}

\newcommand{\e}{\varepsilon}

\newcommand{\polylog}{{\mathrm{polylog}}}

\newcommand{\supp}{{{\mathrm{supp}}}}

\newcommand{\E}{\mathbb{E}}

\DeclareMathOperator{\diam}{diam}

\newtheorem{theorem}{Theorem}[section]
\newtheorem{lemma}[theorem]{Lemma}

\newtheorem{claim}[theorem]{Claim}

\newtheorem{corollary}[theorem]{Corollary}

\newtheorem{definition}[theorem]{Definition}

\newtheorem{fact}[theorem]{Fact}

\newtheorem{remark}{Remark}[section]

\newcommand\card[1]{\left| #1 \right|} 
\begin{document}

\title{Euclidean distortion and the Sparsest Cut}

\author{
 Sanjeev Arora\thanks{Supported by David and Lucile Packard
Fellowship and NSF grant CCR-0205594.  Princeton University.
\tt{arora@cs.princeton.edu}} \and
James R. Lee\thanks{Supported by NSF grant CCR-0121555 and an NSF
Graduate Research Fellowship.  U. C. Berkeley.
\tt{jrl@cs.berkeley.edu}} \and Assaf Naor \thanks{Microsoft
Research.  \tt{anaor@microsoft.com}}}

\date{}
\maketitle

\begin{abstract}
 We prove that every $n$-point metric space of negative
type (and, in particular, every $n$-point subset of $L_1$) embeds
into a Euclidean space with distortion $O(\sqrt{\log n} \cdot\log
\log n)$, a result which is tight up to the iterated logarithm
factor. As a consequence, we obtain the best known polynomial-time
approximation algorithm for the Sparsest Cut problem with general
demands.  Namely, if the demand is supported on a subset of size
$k$, we achieve an approximation ratio of $O(\sqrt{\log k}\cdot \log
\log k)$.
\end{abstract}

\section{Introduction}

Bi-Lipschitz embeddings of finite metric spaces, a topic originally
studied in geometric analysis and Banach space theory, became an
integral part of theoretical computer science following work of
Linial, London, and Rabinovich~\cite{LLR95}. They presented an
algorithmic version of a result of Bourgain~\cite{Bourgain85} which
shows that every $n$-point metric space embeds into $L_2$ with
distortion $O(\log n)$. This geometric viewpoint offers a way to
understand the approximation ratios achieved by linear programming
(LP) and semidefinite programming (SDP) relaxations for cut
problems~\cite{LLR95,AR98}. It soon became apparent that further
progress in understanding SDP relaxations would involve improving
Bourgain's general bound of $O(\log n)$ for $n$-point metric spaces
of negative type. For instance,  the approximation ratio achieved by
a well-known SDP relaxation for the general Sparsest Cut problem is
known to coincide exactly with the best-possible distortion bound
achievable for the embedding of $n$-point metrics of negative type
into $L_1$---a striking connection between pure mathematics and
algorithm design.

Further progress on these problems required new insights into the
structure of metric spaces of negative type, and the design of more
sophisticated and flexible embedding methods for finite metrics.
Coincidentally, significant progress was made recently on both these
fronts. Arora, Rao and Vazirani~\cite{ARV04} proved a new structural
theorem about metric spaces of negative type and used it to design
an $O(\sqrt{\log n})$-approximation algorithm for {\em uniform case}
of the Sparsest Cut problem. Krauthgamer, Lee, Mendel and
Naor~\cite{KLMN04} introduced a new embedding method called {\em
measured descent} which unified and strengthened many existing
embedding techniques, and they used it to solve a number of open
problems in the field.

These breakthroughs indeed resulted in improved embeddings for
negative type metrics; Chawla, Gupta, and R\"acke~\cite{CGR04}
used the structural theorem of~\cite{ARV04} (specifically, its
stronger form due to Lee~\cite{L04}), in conjunction with measured
descent to show that every $n$-point metric of negative type
embeds into $L_2$ with distortion $O(\log n)^{3/4}$. In the
present work, we show how one can achieve distortion $O(\sqrt{\log
n}\cdot \log\log n)$. This almost matches the $35$-year-old lower
bound of $\sqrt{\log n}$ from Enflo~\cite{Enflo69}. Our methods
use the results of~\cite{ARV04,L04,CGR04} essentially as a ``black
box,'' together with an enhancement of the measured descent
technique.

Recall that a metric space $(X,d)$ is said to be of {\em negative
type} if $(X,\sqrt{d})$ is isometric to a subset of Euclidean
space. In particular, it is well known that $L_1$ is of negative
type.  (We also remind the reader that $L_2$ is isometrically
equivalent to a subset of $L_1$.) The parameter $c_2(X)$, known as
the {\em Euclidean distortion} of $X$, is the least distortion
with which $X$ embeds into Hilbert space, i.e. it is the minimum
of $\mathsf{distortion}(f) = \|f\|_{\Lip}\cdot\|f^{-1}\|_{\Lip}$ over
all bijections $f:X\hookrightarrow L_2$. The mathematical
investigation of the problem we study here goes back to the work
of Enflo~\cite{Enflo69}, who showed that the Euclidean distortion
of the Hamming cube $\Omega_d=\{0,1\}^d$ equals
$\sqrt{d}=\sqrt{\log_2|\Omega_d|}$. The following natural question
 is folklore in geometric and functional analysis.

\medskip
{\em
 ``Is the discrete $d$-dimensional hypercube the
most non-Euclidean $2^d$-point subset of $L_1$?"}

\medskip
\noindent A positive answer to this question would imply that any
$n$-point subset of $L_1$ embeds in $L_2$ with distortion
$O(\sqrt{\log n})$. In fact, motivated by F. John's theorem in
convex geometry (see~\cite{MS86}), Johnson and
Lindenstrauss~\cite{jl} asked in 1983 whether every $n$-point metric
space embeds into $L_2$ with distortion $O(\sqrt{\log n})$. Here,
the analogy between finite dimensional normed spaces and finite
metric spaces is not complete: Bourgain~\cite{Bourgain85} has shown
that for any $n$-point metric space $X$, $c_2(X)=O(\log n)$, and
this result is existentially optimal~\cite{LLR95,AR98}. By now we
understand that finite metric spaces (namely expander graphs) can
exhibit an isoperimetric profile which no Euclidean space can
achieve, and this is the reason for the discrepancy with John's
theorem. However, it is known (see~\cite{KLMN04}) that several
natural restricted classes of metrics do adhere to the $O(\sqrt{\log
n})$ Euclidean distortion suggested by John's theorem.  Arguably,
for applications in theoretical computer science, the most important
restricted class of metrics are those of negative type, yet
improvements over Bourgain's theorem for such metrics have long
resisted the attempts of mathematicians and computer scientists.

 The present paper is devoted to proving that up to iterated
logarithmic factors, the answer to the above question is positive.
This yields a general tool for the rounding of certain classes of
semi-definite programs. As a result, we obtain the best-known
polynomial time algorithm for the approximation of the Sparsest
Cut problem with general demands, improving over the previous
bounds due to \cite{CGR04} and the preceding works \cite{LLR95}
and \cite{AR98} (which yield an $O(\log n)$ approximation). This
problem is described in Section~\ref{section:algorithmic}. We now
state our main result.  In the case of metrics of negative
type (and not just $L_1$ metrics), it answers positively (up to
iterated logarithms) a well known conjecture in theoretical
computer science and metric geometry (stated explicitly by Goemans
in~\cite{Goem97}).

\begin{theorem}\label{thm:main}
Let $(X,d)$ be an $n$-point metric space of negative type. Then
$$
c_2(X)=O\left(\sqrt{\log n}\cdot \log\log n\right).
$$
\end{theorem}


\medskip \noindent {\bf Related work.} Until recently, there was
little solid evidence behind the conjecture that any $n$-point subset of
$L_1$ embeds in Hilbert space with distortion $O(\sqrt{\log n})$.
In the paper~\cite{LMN04}, Lee, Mendel and Naor  show that any
$n$-point subset of $L_1$ embeds into Hilbert space with {\em
average distortion} $O(\sqrt{\log n})$. Arora, Rao, and
Vazirani~\cite{ARV04} have shown that $O(\sqrt{\log n})$
distortion is achievable using a different notion of average
distortion, which turns out to be more relevant for bounding the
actual distortion. As described above, combining their result with
the measured descent technique of Krauthgamer, Lee, Mendel and
Naor~\cite{KLMN04}, Chawla, Gupta, and R\"acke~\cite{CGR04} have
recently proved that for any $n$-point metric space $X$ of
negative type, $c_2(X)=O(\log n)^{3/4}$.
It was conjectured \cite[pg. 379]{Mat01} that $n$-point metrics of
negative type embed into $L_1$ with distortion $O(1)$. Recently,
Khot and Vishnoi \cite{KV04} have obtained a lower bound of
$\Omega(\log \log n)^\delta$, for some constant
$\delta>0$.

Our results also suggest that the dimension reduction lower bound
of Brinkman and Charikar~\cite{BC03} (see also~\cite{LNdiamond})
is tight for certain distortions. They show that embedding certain
$n$-point subsets of $L_1$ into $\ell_1^{d}$ with distortion $D$
requires that $d \geq n^{\Omega(1/D^2)}$. Theorem \ref{thm:main},
 together with theorems of Johnson and Lindenstrauss \cite{jl}
 and Figiel, Lindenstrauss, and Milman \cite{flm}, yields
 an embedding of every $n$-point subset of $L_1$ into
 $\ell_1^{O(\log n)}$ with distortion $O\left(\sqrt{\log
n}\cdot\log \log n\right)$.

\subsection{Algorithmic application: The Sparsest Cut problem with general
demands}\label{section:algorithmic}

In this section, we briefly describe an application of
Theorem~\ref{thm:main} to the Sparsest Cut problem with general
demands (and its relation to the {\em multi-commodity flow}
problem). This is a fundamental NP-hard combinatorial optimization
problem---we refer the interested reader to the
articles~\cite{LR99,AKRR90,LLR95,AR98}, the
survey~\cite{Shmoys95}, and Chapter 21 of the book
\cite{Vazirani01} for additional information on Sparsest Cut, and
its applications to the design of approximation algorithms.

Let $G=(V,E)$ be a graph (network), with a {\em capacity} $C(e)\ge
0$ associated to every edge $e\in E$. Assume that we are given $k$
pairs of vertices $(s_1,t_1),...,(s_k,t_k)\in V\times V$ and
$D_1,\ldots,D_k\ge 1$. We think of the $s_i$ as {\em sources}, the
$t_i$ as {\em targets}, and the value $D_i$ as the {\em demand} of
the {\em terminal pair} $(s_i,t_i)$ for some {\em commodity}
$\kappa_i$.
The problem is said to have {\em uniform demands}
if {\em every} pair $u,v \in V$ occurs as some $(s_i, t_i)$ pair
with $D_i = 1$.

In the {\em MaxFlow} problem the objective is to maximize the {\em
fraction} $\lambda$ of the demand that can be shipped
simultaneously for all the commodities, subject to the capacity
constraints. Denote this maximum by $\lambda^*$. A trivial upper
bound on $\lambda^*$ is the {\em cut ratio}.
Given any subset $S \subseteq V$,
we write
$$
\Phi(S) = \frac{\sum_{uv \in E} C(uv) \cdot |{\bf 1}_S(u) - {\bf 1}_S(v)|}
               {\sum_{i=1}^k D_i \cdot |{\bf 1}_S(s_i) - {\bf 1}_S(t_i)|},
$$
where ${\bf 1}_S$ is the characteristic function of $S$.
The value $\Phi^* = \min_{S \subseteq V} \Phi(S)$
is the minimum over all cuts (partitions) of $V$, of the ratio
between the total capacity crossing the cut and the total demand
crossing the cut. In the case of a single commodity (i.e. $k=1$)
the classical MaxFlow-MinCut theorem states that
$\lambda^*=\Phi^*$, but in general this is no longer the case. It
is known~\cite{LR99,LLR95,AR98} that $\Phi^*=O(\log k)\lambda^*$.
This result is perhaps the first striking application of metric
embeddings in combinatorial optimization (specifically, it uses
Bourgain's embedding theorem~\cite{Bourgain85}).

Computing $\Phi^*$ is NP-hard~\cite{SM90}. Moreover, finding a cut
for which $\Phi^*$ is (approximately) attained  is a basic step in
approximation algorithms for several NP-hard
problems~\cite{LR99,AKRR90,Shmoys95}. The best known algorithm for
computing $\Phi^*$ in the case of uniform demands is due
to~\cite{ARV04}, where an approximation ratio of $O(\sqrt{\log n})$
is achieved. In the case of general demands, an approximation ratio
of $O(\log k)^{3/4}$ is obtained in~\cite{CGR04}. Here, as an
application of Theorem~\ref{thm:main}, we prove the following
theorem:

\begin{theorem}\label{thm:algorithm} Using the above notation, there
exists a polynomial-time algorithm which produces a subset
$S\subseteq V$ for which
$$
\Phi(S) = O\left(\sqrt{\log k}\cdot \log\log
k\right)\Phi^*.
$$
\end{theorem}

\bigskip

\noindent{\bf Structure of the paper:} This paper is organized as
follows. In Section~\ref{section:overview} we present an informal
overview of the ideas involved in the proof of
Theorem~\ref{thm:main}. Section~\ref{section:single} is devoted to
various preliminaries on the geometry of metrics of negative type.
Theorem~\ref{thm:main} is proved in Section~\ref{section:proof},
and the algorithm of Theorem~\ref{thm:algorithm} is described and
analyzed in Section~\ref{section:sparsest}. We end with
Section~\ref{section:conclude}, which contains additional remarks
and open problems.

\section{Overview of the proof of Theorem~\ref{thm:main}}\label{section:overview}

\noindent
{\bf Remarks on notation.}
When we write $E \gtrsim F$ for two expressions $E$ and $F$, we intend this to mean that
there exists some $\epsilon > 0$ such that $E \geq \epsilon F$,
where $\epsilon$ is intended to be a universal constant, independent
of the variables or parameters on which $E$ and $F$ depend.

We will often work with Hilbert spaces of the following form: If
$H$ is a Hilbert space, and $(\Omega,\mu)$ is a probability space,
we use $L_2(H, \Omega, \mu)$ to denote the Hilbert space of
$H$-valued random variables $Z$ with norm
$||Z||_{L_2(H,\Omega,\mu)} = \sqrt{\mathbb E\,||Z||_H^2}$. When
$H, \Omega$ are clear from context, we simply write $L_2(\mu)$
and denote $||\cdot||_H$ by $||\cdot||_2$.

\medskip

Our proof of Theorem~\ref{thm:main} has little to do with metrics
of negative type; the connection to such spaces comes through the
techniques of \cite{ARV04,L04,CGR04} and is laid out in
Section~\ref{section:single}. Instead, we present a general
theorem about gluing together various maps from finite metric
spaces into Hilbert spaces (and, more generally, $L_p$ spaces for
$p \in [1,\infty)$).  Our starting point is the following type of
ensemble.

Let $(X,d)$ be an $n$-point metric space.
Suppose that for every $\tau \geq 0$, and every subset $S \subseteq X$,
there exists a 1-Lipschitz map $\f_{S,\tau} : X \to L_2$ with
\begin{equation}
||\f_{S,\tau}(x) - \f_{S,\tau}(y)||_2 \geq \frac{\tau}{\sqrt{\log |S|}}\label{eq:huh}
\end{equation}
whenever $x,y \in S$ and $d(x,y) \in [\tau,2\tau]$.  In general, $\sqrt{\log |S|}$
could be a different function of $|S|$, but we restrict ourselves here for simplicity.
Additionally, let us temporarily define $\f_{\tau} = \f_{X,\tau}$ for every $\tau \geq 0$ so that
for the maps $\{\varphi_\tau\}$, condition \eqref{eq:huh} holds for all $x,y \in X$ and $|S| = n$.
The problem we are now confronted with is how to combine the ensemble of maps $\{\f_{S,\tau}\}$
together to obtain a genuinely bi-Lipschitz map.

There is an obvious approach which comes to mind:  Let $R \subseteq \mathbb Z$
be such that for all $x,y \in X$, there exists $k \in R$ such that
$d(x,y) \in [2^{k}, 2^{k+1}]$.  Now define the map $\varphi : X \to L_2$ by
$\varphi = \bigoplus_{k \in R} \f_{2^k}$.  Clearly we have
both $||\varphi||_{\Lip} \leq \sqrt{R}$ and, for all $x,y \in X$,
$$
||\varphi(x) - \varphi(y)||_2 \geq \frac{d(x,y)}{2 \sqrt{\log n}},
$$
hence $\mathsf{distortion}(\varphi) = O(\sqrt{R \log n})$.  Trivially,
we can choose $R$ so that $|R| \leq n^2$.  A slightly more delicate
argument yields such an $R$ with $|R| \leq O(n)$.  Unfortunately,
we are searching for a bound of the form $\mathsf{distortion}(\varphi) \approx \sqrt{\log n}$,
making this construction useless.

Nevertheless, the key to a better gluing of the given ensemble does lie in the delicate
interplay between the distributions of distances in $X$ and the number of points
in various regions of the space.  The technique of {\em measured descent} from \cite{KLMN04}
relies essentially on two facts about finite metric spaces.  First,
the identity
\begin{equation}\label{eq:ouch}
\sum_{k \in \mathbb Z} \log \frac{|B(x, \alpha\cdot 2^k)|}{|B(x,2^k)|} = O(\log n \log \alpha)
\end{equation}
for any number $\alpha \geq 2$.  (In \cite{KLMN04}, a fixed
constant value of $\alpha$ was used, but for us the quantitative
dependence is crucial, as we will have $\alpha$ depending on $n$.)
This gives a simple bound on the rate that a finite metric space
can expand over all its scales, and is implicitly used in earlier
works under the name of ``region growing'' \cite{LR99,FRT03}.

For the purposes of this description, we will state the second
fact less concretely. Basically, in certain settings, one can
think of the ratio $|B(x, \alpha \tau)|/|B(x,\tau)|$ as the
``local cardinality of the space'' around $x$ at scale $\tau$.  As
an example, if $X = \mathbb R^d$, $B(x,\cdot)$ represents a
Euclidean ball, and $|\cdot|$ is the Lebesgue measure, then this
ratio approximates the number of $\tau$-net points that can be
packed inside a ball of radius $\alpha \tau$. Later, it will
become necessary to randomly partition $X$ into pieces of diameter
at most $2\tau$ while ensuring that pairs $x,y\in X$ with $d(x,y)
\ll \tau$ are usually in the same component of the partition (see
Section \ref{section:definitions} on padded decomposability). It
is known \cite{CKR01,FRT03} that the properties of such partitions
near $x$ depends on the local value $\log \frac{|B(x,
2\tau)|}{|B(x,\tau)|}$.

Following \cite{KLMN04}, this relationship is used in \cite{L04} to prove (roughly) that,
given the maps $\{\varphi_\tau\}_{\tau \geq 0}$ defined above,
there exists a map $\varphi : X \to L_2$ such that
$||\varphi||_\Lip \leq O(\sqrt{\log n})$ and, for $x,y \in X$ with
$d(x,y) \in [2^k, 2^{k+1}]$,
\begin{eqnarray}\label{eq:plumb}
||\varphi(x) - \varphi(y)||_2 &\gtrsim&
\sqrt{\log \frac{|B(x,2^{k+1})|}{|B(x,2^k)|}}
\left(\vphantom{\sqrt{\log \frac{|B(x,2^{k+1})|}{|B(x,2^k)|}}}\right.
\underbrace{\vphantom{\frac{d(x,y)}{\log \frac{|B(x,2^{k+1})|}{|B(x,2^k)|}}}||\varphi_{2^k}(x)-\varphi_{2^k}(y)||}_{(\mathrm{I})}
+ \underbrace{\frac{d(x,y)}{\log \frac{|B(x,2^{k+1})|}{|B(x,2^k)|}}}_{\mathrm{(II)}}
\left.\vphantom{\sqrt{\log \frac{|B(x,2^{k+1})|}{|B(x,2^k)|}}}\right)
\end{eqnarray}
The contribution (II) comes from random partitioning and Rao's technique \cite{Rao99}, and
is valid for any metric space $X$.  Observing that $\mathrm{(I)} \gtrsim d(x,y)/\sqrt{\log n}$,
and using AM-GM in \eqref{eq:plumb}, one arrives at the lower bound
$$||\varphi(x) - \varphi(y)||_2 \gtrsim \frac{d(x,y)}{(\log n)^{\frac{1}{4}}},$$
hence $\mathsf{distortion}(\varphi) \leq O(\log n)^{\frac{3}{4}}$.
While not obvious at present, the identity \eqref{eq:ouch} is what
allows \cite{L04} to get the leading $\sqrt{\cdot}$ factor in \eqref{eq:plumb} while keeping
$||\varphi||_\Lip$ small (see Theorem \ref{thm:descent}).

In order to get the distortion near $O(\sqrt{\log n})$, we have to dispense with
the contribution (II) which is not derived from the ensemble $\{\varphi_{S,\tau}\}$.
Instead, we would like to pass from the ensemble
$\{\varphi_{S,\tau}\}$ to a family of maps $\{\tilde \varphi_\tau : X \to L_2\}$
for which the contribution of (I) in \eqref{eq:plumb} is replaced by
\begin{equation}\label{eq:lie}
||\tilde \varphi_{2^k}(x) - \tilde \varphi_{2^k}(y)||_2 \gtrsim \frac{d(x,y)}{\sqrt{\log \frac{|B(x,2^{k+1})|}{|B(x,2^k)|}}}.
\end{equation}
Clearly this would finish the proof.
Roughly, the construction of $\tilde \varphi_{2^k}$ proceeds as follows.
We first randomly partition $X$ into
components of diameter about $2^k \alpha$ for some appropriately chosen $\alpha = \alpha(n)$.
Writing the random partition as $X = C_1 \cup C_2 \cup \cdots \cup C_m$, we then
derive subsets $\tilde C_i \subseteq C_i$ by randomly sampling
points from each $C_i$.  Then, we use an appropriately constructed
(random) partition of unity to glue the collection of maps $\{\varphi_{\tilde C_i, 2^k}\}_{i = 1}^m$
together.  To ensure that the resulting map still has $||\tilde \varphi_{2^k}||_\Lip \leq O(1)$,
the partition of unity is constructed carefully using properties of the random partition
(this bears some resemblance to the technique of \cite{LN04} for extending Lipschitz functions).

The key to the proof is the way in which the random samples
$\tilde C_i$ are chosen. We have to maintain the property that
$\tilde C_i$ is a ``good representative'' of $C_i$ at scale $2^k$
(i.e. we need that, on average, $\tilde C_i$ is a $2^k$-dense in
$C_i$). On the other hand, we need to maintain the invariant that
if $x \in C_i$, then $$\log |\tilde C_i| \approx \log \frac{|B(x,
\alpha 2^k)|}{|B(x,2^k)|},$$ so that we can achieve a bound
similar to \eqref{eq:lie} (recall that the quality of the map
$\varphi_{S,\tau}$ depends on $|S| = |\tilde C_i|$).
Unfortunately, this is impossible since for distinct $x,x' \in
C_i$, the above ratios can be quite different.   Instead, we have
a number of phases, one for each estimate of the possible ratio
(see the proof of Theorem \ref{thm:main}). For this to work, we
have to give up on achieving \eqref{eq:lie} exactly, and instead
we weave together the {\em inter-scale} (Lemma \ref{lem:loc})
gluing of \eqref{eq:plumb} with the {\em intra-scale} (Theorem
\ref{thm:descent}) gluing of \eqref{eq:lie} to obtain a
nearly-tight bound of $O(\sqrt{\log n}\cdot \log \log n)$.

\section{Single scale embeddings}\label{section:single}

In this section we present Theorem~\ref{thm:arv}, and derive from it
Lemma \ref{lem:single} which is one of the main tools used in the
proof of the Main Theorem (\ref{thm:main}). It is a concatenation of
the result of Arora, Rao, and Vazirani~\cite{ARV04}, its
strengthening by Lee~\cite{L04}, and the ``reweighting" method of
Chawla, Gupta, and R\"acke~\cite{CGR04}, who use it in conjunction
with \cite{KLMN04} to achieve distortion $O\left(\log
n\right)^{\frac{3}{4}}$. For the sake of completeness, we present below a
sketch of the proof of Theorem~\ref{thm:arv}. Complete details can
be found in the full version of~\cite{L04}, where a more general
result is proved; the statement actually holds for metric spaces
which are quasisymmetrically equivalent to subsets of Hilbert space,
and not only for those of negative type. (See \cite{Heinonen01} for
the definition of quasisymmetry; the relevance of such maps to the
techniques of \cite{ARV04} was first pointed out in \cite{NRS04}).
\begin{theorem}\label{thm:arv}
There exist constants $C \geq 1$ and $0 < p < \frac{1}{2}$ such
that for every $n$-point metric space $(Y,d)$ of negative type and
every $\Delta > 0$, the following holds. There exists a
distribution $\mu$ over subsets $U \subseteq Y$ such that for
every $x,y \in Y$ with $d(x,y) \geq \frac{\Delta}{16}$,
$$
\mu \left\{ U : y \in U \textrm{ and } d(x, U) \geq
\frac{\Delta}{C \sqrt{\log n}}\right\} \geq p.
$$
\end{theorem}

\begin{proof}[Proof (sketch)]
Let $g : Y \to \mathbb \ell_2$ be such that $$d(x,y) = ||g(x) -
g(y)||_2^2$$ for all $x,y \in Y$. By~\cite{MN04}, there exists a
map $T : \ell_2 \to \ell_2$ such that $||T(z)||_2 \leq
\sqrt{\Delta}$ for all $z \in \ell_2$ and
$$
\frac{1}{2} \leq
\frac{||T(z)-T(z')||_2}{\min\{\sqrt{\Delta},||z-z'||_2\}} \leq 1.
$$
for all $z,z' \in \ell_2$. As in \cite{L04}, we let $f : Y \to
\mathbb R^n$ be the map given by $f = T \circ g$ (we remark that
this map can be computed efficiently). Then $f$ is a bi-Lipschitz
embedding (with distortion $2$) of the metric space
$\left(Y,\sqrt{\min\{\Delta,d\}}\right)$ into the Euclidean ball of
radius $\Delta$.

Let $0 < \sigma < 1$ be some constant.
The basic idea is to choose a random $u \in S^{n-1}$ 
and define
\begin{eqnarray*}
L_u &=& \left\{ x \in Y : \langle x, u \rangle \leq \tfrac{-\sigma \sqrt{\Delta}}{\sqrt{n}}\right\}, \\
R_u &=& \left\{ x \in Y : \langle x, u \rangle \geq \tfrac{\sigma
\sqrt{\Delta}}{\sqrt{n}}\right\}.
\end{eqnarray*}
One then prunes the sets by iteratively removing any pairs of
nodes $x \in L_u$, $y \in R_u$ with $d(x, y) \leq
\Delta/\sqrt{\log n}$. At the end one is left with two sets $L_u',
R_u'$. The main result of~\cite{ARV04, L04} is that with high
probability (over the choice of $u$), the number of pairs pruned
from $L_u \times R_u$ is not too large.

\remove{
 $L_u'$ and $R_u'$ each
contain $\Omega(n)$ points (In fact, we are using here the
strengthening of this result from~\cite{L04}, where it is shown in
particular that this conclusion holds true for any subset of a
sphere in Hilbert space which is a quasisymmetric image of any
$n$-point metric space- the fact that this is the crucial
geometric condition in~\cite{ARV04} was first proved
in~\cite{NRS04}. In our particular case, the quasisymmetry is very
simple- we are dealing with a bi-Lipschitz image of
$(Y,\sqrt{\min\{\Delta,d\}})$). }

Let $S_{\Delta} = \{ (x,y) \in Y \times Y : d(x,y) \geq
\frac{\Delta}{16} \}.$ The reweighting idea of~\cite{CGR04} is to
apply the above procedure to a weighted version of the point set
as follows. Let $w : Y \times Y \to \mathbb Z^+$ be an
integer-valued weight function on pairs, with $w(x,y) = w(y,x)$,
$w(x,x) = 0$, and $w(x,y) > 0$ only if $(x,y) \in S_{\Delta}$. This
weight function can be viewed as yielding a new set of points
where each point $x$ is replaced by $\sum_{y\in Y} w(x, y)$
copies, with $w(x,y)$ of them corresponding to the pair $(x, y)$.
One could think of applying the above procedure on this new point
set; note that the pruning procedure above may remove some or all
copies of $x$. Then, as observed in~\cite{CGR04}, the theorems of
\cite{ARV04,L04} imply that with high probability, after the
pruning, we still have $$\sum_{x \in L'_u, y \in R'_u} w(x,y) \gtrsim
\sum_{x,y} w(x,y).$$

The distribution $\mu$ mentioned in the statement of the theorem
is defined using a family of $O(\log n)$ weight functions
described below. Sampling from $\mu$ consists of picking a weight
function from this family and a random direction $u \in S^{n-1}$,
and  then forming sets $L_u', R_u'$ as above using the weight
function.
Let us call these sets $L_u'(w), R_u'(w)$.
One then outputs the set $U$ of all points $x$ for which
any ``copy'' falls into $L_u'(w)$.

Now we define the family of weight functions. The initial weight function has
$w_0(x,y) = n^4$ for all $(x,y) \in S_{\Delta}$. Given $w_k$,
obtain $w_{k+1}$ as follows.
If
$$
\mu \left \{ u \in S^{n-1} : (x,y) \in L'_u(w_k) \times R'_u(w_k) \right \} \geq 0.1,
$$
we set $w_{k+1}(x,y) = \frac{1}{2} w_k(x, y)$.
Otherwise, we set $w_{k+1}(x,y) = w_k(x,y)$.
A simple argument (presented in~\cite{CGR04}) shows
that by repeating this $O(\log n)$ times we obtain $O(\log n)$
weight functions such that for every pair $(x, y) \in S_{\Delta}$
the following is true: If one picks a random weight function $w$ and a
random direction $u \in S^{n-1}$, then with constant probability we
have $(x,y) \in L_u'(w) \times R_u'(w)$.
\end{proof}

\subsection{Padded decomposability and random zero
sets}\label{section:definitions}

Theorem~\ref{thm:arv} is the only way the negative type property
will be used in what follows. It is therefore helpful to introduce
it as an abstract property of metric spaces. Let $(X,d)$ be an
$n$-point metric space.

\begin{definition}[Random zero-sets]\label{def:zero} Given $\Delta,\ \zeta>0$,  and $p\in (0,1)$ we say
that $X$ admits a random zero set at scale $\Delta$ which is
$\zeta$-spreading with probability $p$ if there is a distribution
$\mu$ over subsets $Z\subseteq X$ such that for every $x,y\in X$
with $d(x,y)\ge \Delta$,
$$
\mu \left\{Z\subseteq X:\ y\in Z\ \mathrm{ and }\ d(x,Z)\ge
\frac{\Delta}{\zeta}\right\}\ge p.
$$
We denote by $\zeta(X;p)$ the least $\zeta>0$ such that for every
$\Delta>0$, $X$ admits a random zero set at scale $\Delta$ which
is $\zeta$-spreading with probability $p$. Finally, given $k\le n$
we define
$$
\zeta_k(X;p)= \max_{\substack{Y\subseteq X\\|Y|\le k}}\zeta(Y;p).
$$
\end{definition}
With this definition, Theorem~\ref{thm:arv} implies that there
exists a universal constant $p\in (0,1)$ such that for every
$n$-point
 metric space $(X,d)$ of negative type, $\zeta(X;p)=O(\sqrt{\log
 n})$.

\medskip
We now recall the related notion of {\em padded decomposability}.
Given a partition $P$ of $X$ and $x\in X$ we denote by $P(x)\in P$
the unique element of $P$ to which $x$ belongs. In what follows we
sometimes refer to $P(x)$ as the {\em cluster} of $x$.
\begin{definition}[Decomposition bundle, modulus of padded
decomposability]\label{def:bundle} Following~\cite{KLMN04} we say
that $\{P_\Delta\}_{\Delta>0}$ is an $\alpha$-padded decomposition
bundle of a metric space $X$ if for every $\Delta>0$, $P_\Delta$
is a random partition of $X$ (whose distribution we denote by
$\nu$) with the following properties:
\begin{enumerate}
\item For all $P\in \supp(\nu)$ and all $C\in P$ we have that
$\diam(C)< \Delta$. \item For every $x\in X$ we have that
$$\nu\{P:\ B(x,\Delta/\alpha)\subseteq P(x)\}\ge \frac12.$$
\end{enumerate}
The  modulus of padded decomposability of $X$, denoted $\alpha_X$,
is defined as the largest constant $\alpha>0$ such that $X$ admits
an $\alpha$-padded decomposition bundle.
\end{definition}

As observed in~\cite{KLMN04}, the results
of~\cite{LinialSaks,Bartal96} imply that $\alpha_X=O(\log
\card{X})$, and this will be used in the ensuing arguments. The
following useful fact relates the notions of padded decomposability
and random zero sets. Its proof is motivated by an argument of
Rao~\cite{Rao99}.

\begin{fact}\label{fact:connection} $\zeta(X;1/8)\le \alpha_X$.
\end{fact}

\begin{proof} Fix $\Delta>0$ and let $P$ be a partition of $X$
into subsets of diameter less than $\Delta$. Given $x\in X$ we
denote by $\pi_P(x)$ the largest radius $r$ for which
$B(x,r)\subseteq P(x)$. Let $\{\e_C\}_{C\in P}$ be i.i.d.
symmetric $\{0,1\}$-valued Bernoulli random variables. Let $Z_P$
be a random subset of $X$ given by
$$
Z_P=\bigcup_{C\in P:\ \e_C=0} C\ .
$$

If $x,y\in X$ satisfy $d(x,y)\ge \Delta$ then $P(x)\neq P(y)$. It
follows that
$$
\Pr[y\in Z_P\wedge d(x,Z_P)\ge \pi_P(x)]\ge \frac14.
$$
By the definition of $\alpha_X$, there exists a distribution over
partitions $P$ of $X$ into subsets of diameter less than $\Delta$
such that for every $x\in X$ with probability at least $1/2$,
$\pi_P(x)\ge \Delta/\alpha_X$. The required result now follows by
considering the random zero set $Z_P$.
\end{proof}

We end this section with the following simple lemma, which shows
that the existence of random zero sets implies the existence of
embeddings into $L_2$ which are bi-Lipschitz on a fixed distance
scale.

\begin{lemma}[Random zero sets yield single scale
embeddings]\label{lem:single} For every finite metric space $X$,
every $S\subseteq X$ every $p\in (0,1)$ and every $\tau>0$, there
exists a $1$-Lipschitz mapping $\f:X\to L_2$ such that for every
$x,y\in S$ with $d(x,y)\ge \tau$,
$$
\|\f(x)-\f(y)\|_2\ge \frac{\tau\sqrt{p}}{\zeta(S;p)}.
$$
\end{lemma}

\begin{proof}
By the definition of $\zeta(S,p)$ there exists a distribution
$\mu$ over subsets $Z\subseteq S$ such that for every $x,y\in S$
with $d(x,y)\ge \tau$,
$$
\mu \left\{Z\subseteq S:\ y\in Z\ \mathrm{ and }\ d(x,Z)\ge
\frac{\tau}{\zeta(S;p)}\right\}\ge p.
$$
Define $\f:X\to L_2(\mu)$ by $\f(x)=d(x,Z)$. Clearly $\f$ is
$1$-Lipschitz. Moreover, for every $x,y\in S$ with $d(x,y)\ge
\tau$,
\begin{eqnarray*}
\|\f(x)-\f(y)\|_{L_2(\mu)}^2=
\mathbb E_\mu \left[d(x,Z)-d(y,Z)\right]^2 \ge
p\cdot \left(\frac{\tau}{\zeta(S;p)}\right)^2.
\end{eqnarray*}
\end{proof}

\section{Proof of Theorem~\ref{thm:main}}\label{section:proof}

The primary result of this section is the following theorem.

\begin{theorem}\label{thm:simple}
Let $(X,d)$ be an $n$-point metric space.
Suppose there exist constants $C > 0$ and $\frac{1}{2} \leq \varepsilon \leq 1$, such
that for every $\tau \geq 0$, and every subset $S \subseteq X$,
there exists a 1-Lipschitz map $\f_{S,\tau} : X \to L_2$ with
$$
||\f_{S,\tau}(x) - \f_{S,\tau}(y)||_2 \geq \frac{\tau}{C (\log
|S|)^\varepsilon}
$$
whenever $x,y \in S$ and $d(x,y) \in [\tau,6\tau]$.  Then $c_2(X)
\leq O(1)\cdot C(\log n)^{\varepsilon} \log \log n$.
\end{theorem}

Theorem~\ref{thm:simple} implies Theorem~\ref{thm:main}. Indeed,
if $X$ is an $n$ point metric space such that for some $p\in
(0,1)$, $\e\in [1/2,1]$, and $C>0$, we have for every $k\le n$,
$\zeta_k(X;p)\le C(\log k)^\e$, then Theorem~\ref{thm:simple}
together with Lemma~\ref{lem:single} imply that
$$
c_2(X)=O\left(\frac{C(\log n)^\e\log \log n}{\sqrt{p}}\right).
$$
Theorem~\ref{thm:main} follows since by Theorem~\ref{thm:arv} we
know that for some universal constant $p\in (0,1)$, if $X$ is a
metric space of negative type then for all $k$,
$\zeta_k(X;p)=O\left(\sqrt{\log n}\right)$.

\medskip

The proof of Theorem~\ref{thm:simple} will be broken down into
several steps. In what follows we fix a finite metric space $X$,
and for $K \geq 1$, $\tau \geq 0$, define
$$
S_\tau(K) = \left\{ x \in X : \left|B\left(x, 8\tau
\alpha_X\right)\right| \leq K \left|B\left(x, \frac{\tau}{12C(\log
K)^{\varepsilon}}\right)\right|\right\}.
$$

\begin{lemma}[Embedding neighborhoods]
\label{lem:neighborhood}
Let $S \subseteq X$, $\tau \geq 0$, and
assume that there exists a 1-Lipschitz map $\f : X \to L_2$
satisfying
$$
||\f(x) - \f(y)||_2 \geq \frac{\tau}{L}
$$
for $x,y \in S$, $d(x,y) \in[\tau/2,3\tau]$ and some $L \geq 2$.
Then there is a 1-Lipschitz map $h : X \to L_2$ with
$$
||h(x) - h(y)||_2 \geq \frac{\tau}{9L}
$$
whenever $d(x,S) \leq \frac{\tau}{6L}$, $y \in X$, and $d(x,y)
\in[\tau,2\tau]$.
\end{lemma}

\begin{proof}
Define $g : X \to \mathbb R$ by $g(x) = d(x,S)$, and set $h =
\frac{1}{\sqrt{2}}(\f \oplus g)$.  If $d(y, S) > \frac{\tau}{3L}$,
then
$$
||h(x) - h(y)||_2 \geq \frac{1}{\sqrt{2}} ||g(x) - g(y)||_2 \geq
\frac{1}{\sqrt{2}}(d(y,S) - d(x,S)) \geq \frac{1}{\sqrt{2}}\cdot
\frac{\tau}{6L}.
$$
Otherwise, let $x',y' \in S$ be such that $d(x,x') \leq
\frac{\tau}{6L}, d(y,y') \leq \frac{\tau}{3L}$, and observe that
$$d(x',y') \in\left[ d(x,y) - \frac{\tau}{6L} - \frac{\tau}{3L},d(x,y) + \frac{\tau}{6L} + \frac{\tau}{3L}\right] \subseteq \left[\frac{\tau}{2},3\tau\right].$$ Using
our assumptions on $\f$, we have
$$
||\f(x) - \f(y)||_2 \geq ||\f(x') - \f(y')||_2 -
||\f||_{\Lip}\left(\frac{\tau}{6L} + \frac{\tau}{3L}\right) \geq
\frac{\tau}{2L},
$$
hence $||h(x) - h(y)||_2 \geq \frac{1}{\sqrt{2}}\cdot\frac{
\tau}{2L}.$
\end{proof}

\begin{lemma}[Random subsets]
\label{lem:subsets} Assume that $X$ satisfies the conditions of
Theorem~\ref{thm:simple}, and suppose that $U \subseteq X$ and $k
\geq 2$. Define
\begin{equation}\label{eq:Ttau}
T_{\tau}(U;k) = \left\{ x \in U : |U| \leq k\, \left|B\left(x,\frac{\tau}{12C(\log k)^{\varepsilon}}\right)\right| \right\}.
\end{equation}
Then there exists a 1-Lipschitz map $\gamma_{U,k} : X \to L_2$ such that
$$
||\gamma_{U,k}(x) - \gamma_{U,k}(y)||_2 \geq \frac{\tau}{30 C(\log
k)^{\varepsilon}}
$$
whenever $x \in T_{\tau}(U;k), y \in X$ and $d(x,y)\in
[\tau,2\tau]$.
\end{lemma}

\begin{proof}
Let $S$ be a uniformly random subset $S \subseteq U$ with $|S| = \min\{|U|, k\}.$
Let $h_S : X \to L_2$ be the map defined
by $h_S = \frac{1}{\sqrt{2}} (\varphi_{S,\tau/2} \oplus g)$ where $g(x) = d(x,S)$.
Define $\gamma_{U,k} : X \to L_2(L_2, \mu)$,
where $\mu$ is the distribution of the random subset $S$, by
$
\gamma_{U,k}(x) = h_S(x)
$
(recall that $h_S(x)$ is a {\em random} element of $L_2$).  Note
that $\gamma_{U,k}$ is 1-Lipschitz because the same is
true for each $h_S$.

Let $L = 2C (\log |S|)^{\varepsilon}$.
Observe that, by the definition of $T_{\tau}(U;k)$, with
probability at least $1/e$, we have $$S \cap
B\left(x,\frac{\tau}{6L}\right) = S \cap B\left(x, \frac{\tau}{12C
(\log k)^\e}\right) \neq \emptyset.$$ Assuming this holds, we see
that $d(x, S) \leq \frac{\tau}{6L}$. Thus by Lemma
\ref{lem:neighborhood}, $||h_S(x) - h_S(y)||_2 \geq
\frac{\tau}{9L}$. It follows that
$$
||\gamma_{U,k}(x) - \gamma_{U,k}(y)||_2 \geq
\frac{1}{\sqrt{e}}\cdot \frac{\tau}{9L} \geq \frac{\tau}{30C (\log
k)^{\e}}.
$$
\end{proof}

\medskip

In what follows we shall use the fact that for every $\tau>0$
there exists a mapping $G_\tau:L_2\to L_2$ such that for every
$x,y\in L_2$,
\begin{eqnarray}\label{eq:truncation}
\|G_\tau(x)\|_2=\|G_\tau(y)\|_2=\tau\quad \mathrm {and}\quad
 \tfrac12\min\{\tau,\|x-y\|_2\}\le \|G_\tau(x)-G_\tau(y)\|_2\le
\min\{\tau,\|x-y\|_2\}.
\end{eqnarray}
The existence of $G_\tau$ is precisely Lemma 5.2 in~\cite{MN04}.
As in \cite{L04}, we will use the map $G_\tau$ to control
the Lipschitz constant of various functions under partitions of unity.

\begin{lemma}[Localization]\label{lem:loc} Assume that $X$ satisfies the conditions of
Theorem~\ref{thm:simple}. Then for every $\tau \geq 0, k \geq 1$,
there exists a 1-Lipschitz map $\Lambda_{\tau,k} : X \to L_2$ such
that for every $x \in S_\tau(k) ,y \in X$ with $d(x,y) \in
[\tau,3\tau]$,
$$
||\Lambda_{\tau,k}(x) - \Lambda_{\tau,k}(y)||_2 \geq
\frac{\tau}{240C(\log k)^\e}.
$$
\end{lemma}

\begin{proof}
Let $D = 4\tau \alpha_X$ and take $P_D$ to be a random partition
from the $\alpha_X$-padded bundle ensured by
Definition~\ref{def:bundle}. Define a random mapping $\rho : X \to
\mathbb R$ by
$$
\rho(z) = \min\left\{1,\frac{d(z,X\setminus
P_D(z))}{\tau}\right\}.
$$
Clearly $\|\rho\|_{\Lip}\le 1/\tau$. For each $U \in P_D$, let
$\gamma_{U,k}$ be the corresponding map from Lemma
\ref{lem:subsets}. Finally, define a random map
$\Lambda_{\tau,k}:X\to L_2$ by
$$\Lambda_{\tau,k}(z) = \tfrac{1}{2} \rho(z) \cdot \widehat \gamma_{P_D(z),k}(z),$$
where for $f:X\to L_2$ we write $\widehat f = G_{\tau} \circ f$,
where $G_\tau$ is as in~\eqref{eq:truncation}.

We claim that $\|\Lambda_{\tau,k}\|_{\Lip}\le 1$. Indeed, fix $u,v \in X$. If $P_D(u)=P_D(v)=U$ then
\begin{eqnarray*}
||\Lambda_{\tau,k}(u) - \Lambda_{\tau,k}(v)||_2 &\leq&
\tfrac{1}{2} |\rho(u) - \rho(v)| \cdot ||\widehat
\gamma_{P_D(u),k}(u)||_2 +
\tfrac{1}{2} ||\widehat \gamma_{U,k}(u)-\widehat \gamma_{U,k}(v)||_2 \cdot |\rho(v)| \\
&\leq& \tfrac{1}{2} (\tau ||\rho||_{\Lip} + ||\widehat \gamma_{U,k}||_{\Lip})\,d(u,v) \\
&\leq& d(u,v).
\end{eqnarray*}
Otherwise, assume that $P_D(u)\neq P_D(v)$. In particular, $$d(u,v)\ge \max\{d(u,X\setminus P_D(u)),d(v,X\setminus P_D(v))\}.$$
It follows that
\begin{eqnarray*}
||\Lambda_{\tau,k}(u) - \Lambda_{\tau,k}(v)||_2 &\le&
||\Lambda_{\tau,k}(u)||_2 + ||\Lambda_{\tau,k}(v)||_2\\&\le&
\frac{d(u,X\setminus P_D(u))}{2\tau}\cdot \tau +
\frac{d(v,X\setminus P_D(v))}{2\tau}\cdot \tau\\&\le& d(u,v).
\end{eqnarray*}

Now suppose that $x \in S_{\tau}(k), y \in X$, and $d(x,y) \in [\tau,3\tau]$.
Observe that since $\diam(P_D(x)) \leq D$, we have $P_D(x)
\subseteq B(x, 2D)$.  It follows that since $x\in S_\tau(k)$, we have $x\in
T_\tau(P_D(x);k)$ (recall equation \eqref{eq:Ttau}). Moreover, using the defining
property of the $\alpha_X$-padded bundle, with probability at least $\frac12$,
we have $d(x,X \setminus P_D(x))\ge 5\tau$. Since we are assuming that $d(x,y)\le
3\tau$, this implies that $\rho(x)=\rho(y)=1$. It follows that
\begin{eqnarray*}
\mathbb E\,||\Lambda_{\tau,k}(x) - \Lambda_{\tau,k}(y)||_2&\ge&
\tfrac12 \cdot \tfrac{1}{2}\E\, ||\widehat
\gamma_{P_D(x),k}(x)-\widehat \gamma_{P_D(x),k}(y)||_2
\\&\ge& \tfrac1{8} \E\left(\min\left\{|| \gamma_{P_D(x),k}(x)-
\gamma_{P_D(x),k}(y)||_2,\tau\right\}\right)\\&\ge&
\frac{\tau}{240C (\log k)^\e}.
\end{eqnarray*}
Denoting by $(\Omega,\mu)$ the probability space on which
$\Lambda_{\tau,k}$ is defined, we can think of $\Lambda_{\tau,k}$
as a mapping of $X$ into the Hilbert space $L_2(L_2,\mu)$ which
has the required properties.
\end{proof}

The following theorem is a generalization of the Gluing Lemma in \cite{L04}.
In particular, it is important for us that part (2) treats $x$ and $y$ symmetrically,
unlike in \cite{L04}.

\begin{theorem}[Inter-scale gluing]\label{thm:descent}
Given any $n$-point metric space $(X,d)$ and constants $A, B \geq 1$, and
for every $m \in \mathbb Z$, a 1-Lipschitz map $\phi_m : X \to
L_2$, there exists a map $\varphi : X \to L_2$ which satisfies
\begin{enumerate}
\item $||\varphi||_\Lip \leq O(\sqrt{\log n \log (AB)}).$ \item
For every $x,y \in X$
we have
$$||\varphi(x) - \varphi(y)||_2 \geq \max_{m\in \mathbb Z}\left( \sqrt{\left\lfloor \log\frac{|B(x, 2^{m+1} A)|}{|B(x, 2^m/B)|}\right\rfloor} \cdot
\min\left\{\frac{2^m}{B}, ||\phi_{m}(x) -
\phi_{m}(y)||_2\right\}\right).$$
\end{enumerate}
\end{theorem}

\begin{proof}
Let $\rho : X \to \mathbb R_+$ be any $2B$-Lipschitz map with
$\rho \equiv 1$ on $[1/B, 2A]$, and $\rho \equiv 0$ outside
$[1/2B, 4A]$. For $x \in X$ and $t \geq 0$, define
$$
R(x,t) = \sup \{ R : |B(x,R)| \leq 2^t \},
$$
and observe that $R(\cdot,t)$ is 1-Lipschitz for every value of $t$.
And
for each $m \in \mathbb Z$,
define
$$
\rho_{m,t}(x) = \rho\left(\frac{R(x,t)}{2^m}\right).
$$
Write $\widehat \phi_m = G_{2^m/B} \circ \phi_m$, where
$G_{2^m/B}$ is as in~\eqref{eq:truncation}. Now, for each $t \in
\{1,2,\ldots,\lceil \log_2 n\rceil\}$, define $\psi_t : X \to
\ell_2(L_2)$,
$$
\psi_t(x) = \bigoplus_{m \in \mathbb Z} \rho_{m,t}(x) \cdot \widehat \phi_m(x).
$$
Finally, let $\varphi = \psi_1 \oplus \psi_2 \oplus \cdots \oplus \psi_{\lceil \log_2 n \rceil}.$

First, we bound $||\psi_t||_\Lip$ as follows.
\begin{eqnarray*}
||\psi_t(x) - \psi_t(y)||^2_2 = \sum_{\substack{m \in \mathbb Z \\
\rho_{m,t}(x) + \rho_{m,t}(y) > 0}} ||\rho_{m,t}(x) \widehat
\phi_m(x) - \rho_{m,t}(y) \widehat \phi_m(y)||^2_2.
\end{eqnarray*}
The number of non-zero summands above is at most $O(\log A + \log
B)$. Furthermore, each summand can be bounded as follows.
\begin{eqnarray*}
||\rho_{m,t}(x) \widehat \phi_m(x) - \rho_{m,t}(y) \widehat
\phi_m(y)||_2 & \leq &
|\rho_{m,t}(x) - \rho_{m,t}(y)| \cdot ||\widehat \phi_m(x)||_2 + ||\widehat \phi_m(x) - \widehat \phi_m(y)||_2 \cdot |\rho_{m,t}(y)| \\
&\leq & \left(||\rho_{m,t}||_\Lip \cdot \frac{2^m}{B} + ||\phi_m||_{\Lip}\right)\,d(x,y) \\
&\leq& 4\,d(x,y).
\end{eqnarray*}
Thus $||\psi_t||_{\Lip} \leq O(\sqrt{\log(AB)})$.  It follows that $||\varphi||_\Lip \leq O(\sqrt{\log n \log(AB)})$,
as claimed.

It remains to prove the lower bound.  To this end, fix $m\in
\mathbb Z$, $x,y\in X$ and observe that if $\rho_{m,t}(x) = 1$,
then
\begin{eqnarray}\label{eq:contribution1}
\nonumber||\psi_t(x) - \psi_t(y)||_2 &\geq& \|\widehat \phi_m(x) -
\widehat \phi_m(y)||_2-(1-\rho_{m,t}(y))\cdot \|\widehat
\phi_m(y)\|_2\\& \geq& \frac12\min\left\{\frac{2^m}{B},
||\phi_{m}(x) - \phi_{m}(y)||_2\right\}-\frac{2^m}{B}\cdot
(1-\rho_{m,t}(y)).
\end{eqnarray}
On the other hand
\begin{eqnarray}\label{eq:contribution2}
||\psi_t(x) - \psi_t(y)||_2\ge
||\psi_t(x)||_2-||\psi_t(y)||_2=\|\widehat
\phi_m(x)\|_2-\rho_{m,t}(y)\|\widehat
\phi_m(y)\|_2=\frac{2^m}{B}\cdot (1-\rho_{m,t}(y)).
\end{eqnarray}
Averaging~\eqref{eq:contribution1} and~\eqref{eq:contribution2} we
get that
\begin{eqnarray}\label{eq:contribution}
||\psi_t(x) - \psi_t(y)||_2\ge \frac14\min\left\{\frac{2^m}{B},
||\phi_{m}(x) - \phi_{m}(y)||_2\right\}.
\end{eqnarray}
Hence it suffices to count the number of values of $t$ for which
$\rho_{m,t}(x) = 1$. By our definitions we have that
$$
\rho_{m,t}(x) = 1 \iff \frac{2^m}{B} \leq R(x,t) \leq 2^{m+1} A
\iff t \in [\log |B(x,2^m/B)|,\log |B(x,2^{m+1} A)|].
$$
This completes the proof since the lower
bound~\eqref{eq:contribution} holds for $\left\lfloor \log
\frac{|B(x,2^{m+1} A)|}{|B(x,2^m/B)|} \right\rfloor$ values of
$t$.
\end{proof}

We also present the following base case.

\begin{claim}[Small ratios]\label{claim:small}
Let $X$ be an $n$-point metric space and $\tau,\lambda \geq 0$.
Define the subset
$$
S_{\lambda} = \{ x \in X : |B(x,\tau)| \leq \lambda |B(x, \tau/2)| \}.
$$
Then there exists a 1-Lipschitz map $F : X \to L_2$ such that if
$x \in S_{\lambda}$ and $y \in X$ with $d(x,y) \geq \tau$, then
$$
||F(x) - F(y)|| \geq \frac{\epsilon(\lambda) \tau}{\sqrt{\log n}}.
$$
where $\epsilon(\lambda) > 0$ is a constant depending only on $\lambda$.
\end{claim}

\begin{proof}
For each $t \in \{1,2,\ldots,\lceil \log n \rceil \}$, let $W_t
\subseteq X$ be a random subset which contains each point of $X$
independently with probability $2^{-t}$. Let $g_t(x) =
\min\{d(x,W_t),\tau/4\}$ and define the {\em random} map $f =
\frac{1}{\sqrt{\lceil \log n\rceil}} \left(g_1 \oplus \cdots
\oplus g_{\lceil \log n \rceil}\right)$ so that $||f||_\Lip \leq
1$.  Finally, we define $F : X \to L_2(\mu)$ by $F(x) = f(x)$,
where $\mu$ is the distribution over which the random subsets
$\{W_t\}$ are defined.

Now fix $x \in S_\lambda$ and let $t \in \mathbb N$ be such that
$2^t \leq |B(x,\tau/2)| \leq 2^{t+1}$. Let $\mathcal
E_{\mathrm{far}}$ be the event $\left\{d(x,W_t) \geq
\tau/4\right\}$ and let $\mathcal E_{\mathrm{close}}$ be the event
$\left\{d(x, W_t) \leq \tau/8\right\}$. Clearly both such events
are independent of the values $\{ g_t(z) : d(x,z) \geq \tau \}$
(this relies crucially on the use of $\min \{ \cdot, \tau/4 \}$ in
the definition of $g_t$). In particular, these events are
independent of the value $g_t(y)$. It follows that
\begin{eqnarray*}
\| F(x) - F(y) \|^2_{L_2(\mu)}
&=& \mathbb E_{\mu} \left\|f(x) - f(y)\right\|_2^2 \\
&\gtrsim &\frac{1}{\log n} \mathbb E_{\mu} \left|g_t(x) - g_t(y)\right|_2^2 \\
&\gtrsim& \frac{\tau^2}{\log n} \cdot {\min \big\{ \Pr(\mathcal
E_{\mathrm{far}}), \Pr(\mathcal E_{\mathrm{close}}) \big\}}.
\end{eqnarray*}
Finally, we observe that $\Pr(\mathcal E_{\mathrm{far}})$ and
$\Pr(\mathcal E_{\mathrm{close}})$ can clearly be lower bounded by some $\epsilon(\lambda) > 0$.
\remove{
Clearly the events
$$
\{ \mathcal E_y : d(y,x) \geq R\} \textrm{ and } \{ \mathcal E_y : d(y,x) < R\}
$$
are independent.
}
\end{proof}

\medskip

We are now in position to conclude the proof of
Theorem~\ref{thm:simple}

\begin{proof}[Proof of Theorem~\ref{thm:simple}]
We claim that for every $K \in [2,n]$ there exists a map $f_K : X
\to L_2$ which satisfies
\begin{enumerate}
\item $||f_K||_{\Lip} \leq O(\sqrt{\log n\cdot \log\log n})$.
\item For every $m\in \mathbb Z$ and $x \in S_{2^m}(K), y \in X$
we have
$$
||f_K(x) - f_K(y)||^2_2 \gtrsim \left\lfloor
\log\frac{|B(x, 2^{m+3} \alpha_X)|}{|B(x, 2^m/[12C(\log
K)^{\varepsilon}])|}\right\rfloor\cdot \frac{2^{2m}}{C^2(\log
K)^{2\e }}.
$$
\end{enumerate}
Indeed, $f_K$ is obtained from an application of
Theorem~\ref{thm:descent} to the mappings
$\{\Lambda_{2^m,K}\}_{m\in \mathbb Z}$ from Lemma~\ref{lem:loc}
with $A=4\alpha_X$ and $B=12C(\log K)^{\varepsilon}$ (and using
the fact that $\alpha_X=O(\log n)$).

Observe that for every $m\in \mathbb Z$, $S_{2^m}(n)=X$. Hence,
defining $K_0=n$ and $K_{j+1}=\sqrt{K_j}$, as long as $K_j\ge 4$,
we obtain mappings $f_0,\ldots,f_j:X\to L_2$ satisfying
\begin{enumerate}
\item $\|f_j\|_{\Lip} \leq O(\sqrt{\log n\cdot \log \log n})$.
\item For all $x \in S_{2^m}(K_j)\setminus S_{2^m}(K_{j+1})$ and
$y \in X$ such that $d(x,y) \in [2^{m}, 2^{m+1}]$ and
we have
\begin{eqnarray}
||f_{j}(x) - f_{j}(y)||^2_2 &\gtrsim& \left\lfloor
\log\frac{|B(x, 2^{m+3} \alpha_X)|} {|B(x, 2^m/[12C(\log
K_j)^{\varepsilon}])|}\right\rfloor\cdot \frac{2^{2m}}{C^2(\log
K_j)^{2\e}}\nonumber\\\label{eq:monotonicity} &\gtrsim&
\left\lfloor \log\frac{|B(x, 2^{m+3} \alpha_X)|}
{|B(x, 2^m/[12C(\log K_{j+1})^{\varepsilon}])|}\right\rfloor\cdot
\frac{d(x,y)^2}{C^2(\log K_j)^{2\e}}\\\label{eq:not in
Kj+1} &\gtrsim& \left\lfloor \log K_{j+1}\right\rfloor\cdot
\frac{d(x,y)^2}{C^2(\log
K_j)^{2\e}}\\
&\gtrsim& \frac{d(x,y)^2}{C^2(\log
K_j)^{2\e-1}}\label{eq:more than 4},
\end{eqnarray}
\end{enumerate}
where in~\eqref{eq:monotonicity} we used the fact that $K_{j+1}\le
K_j$ and $d(x,y)\le 2^{m+1}$, in~\eqref{eq:not in Kj+1} we used
the fact that $x\notin S_{2^{m}}(K_{j+1})$, and in~\eqref{eq:more
than 4} we used the fact that $K_{j+1}=\sqrt{K_j}\ge 2$.

This procedure ends after $N$ steps, where $N\le O(\log \log n)$.
Every $x\in S_{2^m}(K_N)$ satisfies
$$
|B(x,2^{m+3}\alpha_X)|\le 4|B(x,2^{m}/[12C])|.
$$
By Claim \ref{claim:small}, there is a mapping $f_{N+1}: X\to
L_2$ which is Lipschitz with constant $O(\sqrt{\log n})$ and for
every $x,y\in S_{2^m}(K_N)$, $\|f_{N+1}(x)-f_{N+1}(y)\|_2\gtrsim
d(x,y)$.

 Consider the map $\Phi=\bigoplus_{j=0}^{N+1}f_j$, which is
Lipschitz with constant $O\left(\sqrt{\log n}\cdot\log \log
n\right)$. For every $x,y\in X$ choose $m\in \mathbb Z$ such that
$d(x,y)\in [2^{m},2^{m+1}]$. If $x,y\in S_{2^m}(K_N)$ then
$$
\|\Phi(x)-\Phi(y)\|_2\ge \|f_{N+1}(x)-f_{N+1}(y)\|_2\gtrsim
d(x,y).
$$
Otherwise, without loss of generality there is $j\in \{0,\ldots,
N-1\}$  such that $x\in S_{2^m}(K_j)\setminus S_{2^m}(K_{j+1})$,
in which case by~\eqref{eq:more than 4}
$$
\|\Phi(x)-\Phi(y)\|_2\ge \|f_{j+1}(x)-f_{j+1}(y)\|_2 \gtrsim
\frac{d(x,y)}{C(\log K_j)^{\e-\frac12}}\ge
\frac{d(x,y)}{C(\log n)^{\e-\frac12}}.
$$
\end{proof}

\section{The sparsest cut problem with general demands}
\label{section:sparsest}

This section is devoted to the proof of
Theorem~\ref{thm:algorithm}. Our argument follows the well known
approach for deducing the algorithmic Theorem~\ref{thm:algorithm}
from the embedding result contained in Theorem~\ref{thm:main} (see
e.g.~\cite{LLR95,AR98,Goem97}).

\subsection{Computing the Euclidean distortion}

In this section, we remark that the maps used to prove Theorem
\ref{thm:main} have a certain ``auto-extendability'' property
which will be used in the next section.  We also recall that it is
possible to find near-optimal Euclidean embeddings using
semi-definite programming~\cite{LLR95}.

\begin{corollary}\label{cor:distortion}
Let $(Y, d)$ be an arbitrary metric space, and fix a $k$-point
subset $X \subseteq Y$. If the space $(X,d)$ is a metric of
negative type, then there exists a 1-Lipschitz map $f : Y \to L_2$
such that the map $f|_X : X \to L_2$ has distortion
$O\left(\sqrt{\log k}\cdot \log \log k\right)$.
\end{corollary}

\begin{proof}
We observe that the maps used to prove Theorem \ref{thm:main},
i.e. those produced in Lemma \ref{lem:single} and Claim \ref{claim:small}, are of
Fr\'echet-type.  In other words, there is a
probability space $(\Omega, \mu)$ over subsets $A_\omega \subseteq
X$ for $\omega \in \Omega$, and we obtain a maps $\varphi_{S,\tau} : X \to
L_2(\mu)$ given by $\varphi_{S,\tau}(x)(\omega) = d(x, A_\omega)$.  We can then
define the extension $\varphi_{S,\tau} : Y \to L_2(\mu)$ by
$$\varphi_{S,\tau}(y)(\omega) = d(y, A_\omega).$$
Thus by extending the ensemble of maps $\{\varphi_{S,\tau}\}$
to the larger space $Y$
before the application of Theorem \ref{thm:simple},
we can ensure that the final embedding
is 1-Lipschitz on $Y$.
\end{proof}

\remove{\begin{theorem}
Let $(Y, d)$ be an arbitrary metric space, and fix a $k$-point
subset $X \subseteq Y$. If the space $(X,d)$ is a metric of
negative type, then there exists a probability space $(\Omega,
\mu)$, and a map $f : Y \to L_2(\Omega, \mu)$ such that
\begin{enumerate}
\item For every $\omega \in \Omega$, the map $x \mapsto
f(x)(\omega)$ is 1-Lipschitz. \item For every $x,y \in X$,
$$
||f(x) - f(y)||_2 \geq \frac{d(x,y)}{C \sqrt{\log k} \log \log k},
$$
\end{enumerate}
where $C > 0$ is some universal constant.
\end{theorem}

\begin{proof}
We observe that the maps used to prove Theorem \ref{thm:main},
i.e. those produced in Lemma \ref{lem:single} and Claim \ref{claim:base}), are of
Fr\'echet-type.  In other words, there is a
probability space $(\Omega, \mu)$ over subsets $A_\omega \subseteq
X$ for $\omega \in \Omega$, and we obtain a map $g : X \to
L_2(\mu)$ given by $g(x)(\omega) = d(x, A_\omega)$. We can then
define the extension $f : Y \to L_2(\mu)$ by
$$f(y)(\omega) = d(y, A_\omega).$$
This ensures that the map $x \mapsto f(x)(\omega)$ is 1-Lipschitz
on $Y$ for every $\omega \in \Omega$.
\end{proof}
}


Now we suppose that $(Y,d)$ is an $n$-point metric space and $X
\subseteq Y$ is a $k$-point subset.

\begin{claim}
There exists a polynomial-time algorithm (in terms of $n$) which,
given $X$ and $Y$, computes a map $f : Y \to L_2$ such that $f|_X$
has minimal distortion among all 1-Lipschitz maps $f$.
\end{claim}

\begin{proof}
We give a semi-definite program computing the optimal $f$, which
can be solved in polynomial time using the methods
of~\cite{GLS81}.

\bigskip

\fbox{
\begin{minipage}{0.9\columnwidth}
$$
\begin{array}{rll}
 & \underline{\mathbf{SDP\ (5.1)}}  \ \ \ \ \ \ \ \ \ \ \ \ \  & \\
 \\
\textrm{max}\qquad & \varepsilon  &  \\
\textrm{s.t.}\qquad & x_u \in \mathbb R^n \qquad\qquad\qquad\qquad & \forall u \in Y \\
                    & \|x_u - x_v\|_2^2 \leq d(u,v)^2 & \forall u,v \in Y \\
                    & \|x_u - x_v\|_2^2 \geq \varepsilon\, d(u,v)^2 & \forall u,v \in X
                    \\
\end{array}
$$
\end{minipage}
}

\end{proof}

\subsection{The Sparsest Cut}

Let $V$ be an $n$-point set with two symmetric weights on pairs
$w_N, w_D : V \times V \to \mathbb R_+$ (i.e. $w_N(x,y) =
w_N(y,x)$ and $w_D(x,y) = w_D(y,x)$). For a subset $S \subseteq
V$, we define the {\em sparsity of $S$} by
$$
\Phi_{w_N, w_D}(S) = \frac{\sum_{u \in S, v \in V\setminus  S}
w_N(u,v)}{\sum_{u \in S, v \in V\setminus S} w_D(u,v)},
$$
and we let $\Phi^*(V,w_N,w_D) = \min_{S \subseteq V} \Phi_{w_N,
w_D}(S)$. (The set $V$ is usually thought of as the vertex set of
a graph with $w_N(u,v)$ supported only on edges $(u,v)$, but this
is unnecessary since we allow arbitrary weight functions.)

Computing the value of $\Phi^*(V,w_N,w_D)$ is NP-hard~\cite{SM90}.
The following semi-definite program is well known to be a
relaxation of $\Phi^*(V,w_N,w_D)$ (see e.g.~\cite{Goem97}).

\bigskip

\fbox{
\begin{minipage}{0.9\columnwidth}
$$
\begin{array}{rll}
 & \underline{\mathbf{SDP\ (5.2)}} & \\
\\
\min &  \sum_{u,v \in V} w_N(u,v) \,\|x_u - x_v\|_2^2  & \\
\textrm{s.t.}
     & x_u \in \mathbb R^n \quad \quad \quad \quad \forall u \in V \\
     & \sum_{u,v \in V} w_D(u,v) \, \|x_u - x_v\|_2^2 = 1 &  \\
     & \|x_u - x_v\|_2^2 \leq \|x_u - x_w\|_2^2 + \|x_w - x_v\|_2^2 \\
     & \qquad\qquad\qquad\quad\,\,\,\forall u,v,w \in V
\end{array}
$$
\end{minipage}}
\bigskip

\noindent Furthermore, an optimal solution to this SDP can be
computed in polynomial time~\cite{GLS81,Goem97}.

\medskip
\noindent {\bf The algorithm.} We now give our algorithm for
rounding SDP (5.2). Suppose that the weight function $w_D$ is
supported only on pairs $u,v$ for which $u,v \in U \subseteq V$,
and let $k = |U|$. Denote $M = 20 \log n$.
\begin{enumerate}
\setlength{\itemsep}{0pt}\setlength{\parsep}{0pt}
\setlength{\topsep}{0pt}\setlength{\partopsep}{0pt} \item Solve
SDP (5.2), yielding a solution $\{x_u\}_{u \in V}\subseteq
\mathbb{R}^{n}$. \item Consider the metric space $(V,d)$ given by
$d(u,v) = \|x_u - x_v\|_2^2$. \item Applying SDP (5.1) to $U$ and
$(V,d)$ (where $Y = V$ and $X = U$), compute the optimal map $f :
V \to \mathbb R^n$. \item Choose $\beta_1, \ldots, \beta_M \in
\{-1,+1\}^n$ independently and uniformly at random. \item For each
$1 \leq i \leq M$, arrange the points of $V$ as $v_1^i, \ldots,
v_n^i$ so that
$$\langle \beta_i, f(v_j^i)\rangle \leq \langle \beta_i, f(v_{j+1}^i)\rangle \textrm{ for
        each $1 \leq j \leq n-1$.}$$
\item Output the sparsest of the $Mn$ cuts
\begin{eqnarray*}
(\{v_1^i, \ldots, v_m^i\}, \{v_{m+1}^i, \ldots, v_n^i\}), \qquad 1
\leq m \leq n-1,\  1 \leq i \leq M.
\end{eqnarray*}
\end{enumerate}

\begin{claim}
With constant probability over the choice of $\beta_1, \ldots,
\beta_M$, the cut $(S, V\setminus S)$ returned by the algorithm
has
\begin{equation}
\Phi_{w_N,w_D}(S) \leq O\left(\sqrt{\log k} \log \log k\right)\,
\Phi^*(V,w_N,w_D).\label{eq:sparse}
\end{equation}
\end{claim}

\begin{proof}
Let $S \subseteq \mathbb R^n$ be the image of $V$ under the map
$f$. Consider the map $g : S \to \mathbb \ell_1^M$ given by $g(x)
= (\langle \beta_1, x \rangle, \ldots, \langle \beta_M, x
\rangle)$. It is well-known (see, e.g. \cite{Achlioptas03, MS86})
that, with constant probability over the choice of
$\{\beta_i\}_{i=1}^M \subseteq S^{n-1}$, $g$ has distortion $O(1)$
(where $S$ is equipped with the Euclidean metric). In this case,
we claim that \eqref{eq:sparse} holds.

To see this, let $S_1, S_2, \ldots, S_{Mn} \subseteq V$ be the
$Mn$ cuts which are tested in line (6).  It is a standard
fact~\cite{LLR95,DL97} that there exist constants $\alpha_1,
\alpha_2, \ldots, \alpha_{Mn}\ge 0$ such that for every $x,y \in
V$,
$$||g(f(x)) - g(f(y))||_1 = \sum_{i=1}^{Mn} \alpha_i
\rho_{S_i}(x,y),$$ where $\rho_{S_i}(x,y) = 1$ if $x$ and $y$ are
on opposite sides of the cut $(S_i, V\setminus S_i)$ and
$\rho_{S_i}(x,y) = 0$ otherwise.

Assume (by scaling) that $g \circ f : Y \to \ell_1^M$ is
1-Lipschitz. Let $\Lambda$ be the distortion of $g \circ f$. By
Corollary \ref{cor:distortion}, $\Lambda = O\left(\sqrt{\log k}
\log \log k\right)$. Recalling that $w_D(u,v) > 0$ only when $u,v
\in U$,
\begin{eqnarray*}
\Phi^*(V,w_N,w_D) &\geq& \frac{\sum_{u,v \in V} w_N(u,v) \|x_u - x_v\|_2^2}{\sum_{u,v \in U} w_D(u,v) \|x_u-x_v\|_2^2} \\
&\geq & \frac{1}{\Lambda} \frac{\sum_{u,v \in V} w_N(u,v)\,
||g(f(u)) - g(f(v))||_1}{\sum_{u,v \in U} w_D(u,v)
\,||g(f(u)) - g(f(v))||_1} \\
&= & \frac{1}{\Lambda} \frac{\sum_{i=1}^{Mn} \alpha_i \sum_{u,v
\in V} w_N(u,v) \rho_{S_i}(u,v)}{\sum_{i=1}^{Mn} \alpha_i
\sum_{u,v \in U} w_D(u,v)
\rho_{S_i}(u,v)} \\
&\geq & \frac{1}{\Lambda} \min_{i} \frac{\sum_{u,v \in V} w_N(u,v)
\rho_{S_i}(u,v)}{\sum_{u,v \in U} w_D(u,v) \rho_{S_i}(u,v)}\\& =&
\frac{\Phi_{w_N,w_D}(S)}{\Lambda}.
\end{eqnarray*}
This completes the proof.
\end{proof}

\section{Concluding remarks}
\label{section:conclude}

\begin{itemize}
\item There are two factors of $O\left(\sqrt{\log \log n}\right)$ which keep
our bound from being optimal up to a constant factor. One factor
of $\sqrt{\log \log n}$ arises because Theorem \ref{thm:descent}
is applied with $A,B \sim \mathrm{polylog}(n)$.
The need for such values  arises out of a certain
non-locality property which seems inherent to the method of proof
in \cite{ARV04}.  We remark that achieving $A = O(1)$ is probably possible,
and it seems that $B$ is the difficult factor.

The other factor arises because, in proving Theorem
\ref{thm:main}, we invoke Theorem \ref{thm:descent} for
$O(\log \log n)$ different values of the parameter $K$.  It is
likely removable by a more technical induction, but we chose to
present the simpler proof.

\item It is an interesting open problem  to understand the exact
distortion required to embed $n$-point negative type metrics into
$L_1$. As mentioned before, the best known lower bound is
$\Omega(\log \log n)^{\delta}$ \cite{KV04}. We also note that
assuming a strong form of the {\em Unique Games Conjecture} is true,
the general Sparsest Cut problem is hard to approximate within a
factor of $\Omega\left(\log \log n\right)$ \cite{SKKRS, KV04}.

\item For the uniform case of Sparsest Cut, it is possible to
achieve a $O\left(\sqrt{\log n}\right)$ approximation in quadratic
time without solving an SDP~\cite{AHK04}. Whether such an algorithm
exists for the general case is an open problem.

\item There is no asymptotic advantage in embedding $n$-point
negative type metrics into $L_p$ for some $p\in (1,\infty)$, $p\neq
2$ (observe that since $L_2$ is isometric to a subset of $L_p$ for
all $p\ge 1$, our embedding into Hilbert space is automatically also
an embedding into $L_p$). Indeed, for $1<p<2$ it is shown
in~\cite{LNdiamond} that there are arbitrarily large $n$-point
subsets of $L_1$ that require distortion
$\Omega\left(\sqrt{(p-1)\log n}\right)$ in any embedding into $L_p$.
For $2<p<\infty$ it follows from~\cite{NS02,NPSS04} that there are
arbitrarily large $n$-point subsets of $L_1$ whose minimal
distortion into $L_p$ is $1+\Theta\left(\sqrt{\frac{\log
n}{p}}\right)$ (the dependence on $n$ follows from~\cite{NS02}, and
the optimal dependence on $p$ follows from the results
of~\cite{NPSS04}). Thus, up to multiplicative constants depending on
$p$ (and the double logarithmic factor in Theorem~\ref{thm:main}),
our result is optimal for all $p\in (1,\infty)$.

\item Let $(X,d_X), (Y,d_Y)$ be metric spaces and
$\eta:[0,\infty)\to [0,\infty)$ a strictly increasing function. A
one to one mapping $f:X\hookrightarrow Y$ is called a {\em
quasisymmetric embedding with modulus~$\eta$} if for every
$x,a,b\in X$ such that $x\neq b$,
$$
\frac{d_Y(f(x),f(a))}{d_Y(f(x),f(b))}\le
\eta\left(\frac{d_X(x,a)}{d_X(x,b)}\right).
$$
We refer to~\cite{Heinonen01} for an account of the theory of
quasisymmetric embeddings. Observe that metrics of negative type
embed quasisymmetrically into Hilbert space. It turns out that our
embedding result generalizes to any $n$ point metric space which
embeds quasisymmetrically into Hilbert space. Indeed, if $(X,d)$
embeds quasisymmetrically into $L_2$ with modulus $\eta$ then, as
shown in the full version of~\cite{L04}, there exists constants
$p=p(\eta)$ and $C=C(\eta)$, depending only on $\eta$, such that
$\zeta(X;p)\le C\sqrt{\log n}$.

\remove{ \item It is worthwhile to note here that if a Banach
space $X$ embeds uniformly into Hilbert space (i.e. there exists a
homeomorphism of $X$ into Hilbert space that is uniformly
continuous in both directions), then $X$ embeds quasisymmetrically
into Hilbert space (and thus our main embedding results holds for
$n$-point subsets of $X$). } \remove{Indeed, by a theorem of
Aharoni, Maurey and Mityagin~\cite{AMM85}, $X$ is linearly
isomorphic to a subspace of $L_0$. By a theorem of
Niki\v{s}hin~\cite{Nikisin72} it follows that $X$ is linearly
isomorphic to a subset of $L_{1/2}$, i.e. there exists a linear
operator $T:X\hookrightarrow L_{1/2}$ and a constant $C>0$ such
that for every $x\in X$, $\|x\|\le \|Tx\|_{1/2}\le C\|x\|$. It is
easy to check (see e.g.~\cite{MN04}) that the metric space
$(L_{1/2},\|x-y\|_{1/2}^{1/2})$ is isometric to a subset of $L_1$.
Thus $(X,\|x-y\|^{1/4})$ is bi-Lipschitz equivalent to a subset of
Hilbert space, as required. }

\end{itemize}

\section*{Acknowledgements}

We thank Satish Rao for useful discussions.

\bibliographystyle{abbrv}

\begin{thebibliography}{10}

\bibitem{Achlioptas03}
D.~Achlioptas.
\newblock Database-friendly random projections: {J}ohnson-{L}indenstrauss with
  binary coins.
\newblock {\em J. Comput. System Sci.}, 66(4):671--687, 2003.
\newblock Special issue on PODS 2001 (Santa Barbara, CA).

\bibitem{AKRR90}
A.~Agrawal, P.~Klein, R.~Ravi, and S.~Rao.
\newblock Approximation through multicommodity flow.
\newblock In {\em 31st Annual Symposium on Foundations of Computer Science},
  pages 726--737. IEEE Computer Soc., Los Alamitos, CA, 1990.

\bibitem{AHK04}
S.~Arora, E.~Hazan, and S.~Kale.
\newblock ${O}(\sqrt{\log n})$ approximation to {SPARSEST CUT} in
  $\tilde{O}(n^2)$ time.
\newblock In {\em 45th Annual Syposium on Foundations of Computer Science},
  pages 238--247. IEEE Computer Society, 2004.

\bibitem{ARV04}
S.~Arora, S.~Rao, and U.~Vazirani.
\newblock Expander flows, geometric embeddings, and graph partitionings.
\newblock In {\em 36th Annual Symposium on the Theory of Computing}, pages
  222--231, 2004.

\bibitem{AR98}
Y.~Aumann and Y.~Rabani.
\newblock An {$O(\log k)$} approximate min-cut max-flow theorem and
  approximation algorithm.
\newblock {\em SIAM J. Comput.}, 27(1):291--301 (electronic), 1998.

\bibitem{Bartal96}
Y.~Bartal.
\newblock Probabilistic approximations of metric space and its algorithmic
  application.
\newblock In {\em 37th Annual Symposium on Foundations of Computer Science},
  pages 183--193, Oct. 1996.

\bibitem{Bourgain85}
J.~Bourgain.
\newblock On {L}ipschitz embedding of finite metric spaces in {H}ilbert space.
\newblock {\em Israel J. Math.}, 52(1-2):46--52, 1985.

\bibitem{BC03}
B.~Brinkman and M.~Charikar.
\newblock On the impossibility of dimension reduction in $\ell_1$.
\newblock In {\em Proceedings of the 44th Annual IEEE Conference on Foundations
  of Computer Science}, pages 514--523, 2003.

\bibitem{CKR01}
G.~Calinescu, H.~Karloff, and Y.~Rabani.
\newblock Approximation algorithms for the 0-extension problem.
\newblock In {\em Proceedings of the 12th Annual ACM-SIAM Symposium on Discrete
  Algorithms}, pages 8--16, Philadelphia, PA, 2001.

\bibitem{CGR04}
S.~Chawla, A.~Gupta, and H.~R\"acke.
\newblock Embeddings of negative-type metrics and an improved approximation to
  generalized sparsest cut.
\newblock In {\em Proceedings of the 16th Annual ACM-SIAM Symposium on Discrete
  Algorithms}, pages 102--111, Vancouver, 2005.

\bibitem{SKKRS}
S.~Chawla, R.~Krauthgamer, R.~Kumar, Y.~Rabani, and D.~Sivakumar.
\newblock On embeddability of negative type metrics into $\ell_1$.
\newblock In {\em Proceedings of the 20th Annual IEEE Conference on
  Computational Complexity}, pages 144--153, 2005.

\bibitem{DL97}
M.~M. Deza and M.~Laurent.
\newblock {\em Geometry of cuts and metrics}, volume~15 of {\em Algorithms and
  Combinatorics}.
\newblock Springer-Verlag, Berlin, 1997.

\bibitem{Enflo69}
P.~Enflo.
\newblock On the nonexistence of uniform homeomorphisms between
  ${L}\sb{p}$-spaces.
\newblock {\em Ark. Mat.}, 8:103--105, 1969.

\bibitem{FRT03}
J.~Fakcharoenphol, S.~Rao, and K.~Talwar.
\newblock A tight bound on approximating arbitrary metrics by tree metrics.
\newblock In {\em Proceedings of the 35th Annual ACM Symposium on Theory of
  Computing}, pages 448--455, 2003.

\bibitem{flm}
T.~Figiel, J.~Lindenstrauss, and V.~D. Milman.
\newblock The dimension of almost spherical sections of convex bodies.
\newblock {\em Acta Math.}, 139(1-2):53--94, 1977.

\bibitem{Goem97}
M.~X. Goemans.
\newblock Semidefinite programming and combinatorial optimization.
\newblock {\em Mathematical {P}rogramming}, 49:143--161, 1997.

\bibitem{GLS81}
M.~Gr\"otschel, L.~L.~Lov\'asz, and A.~Schrijver.
\newblock The ellipsoid method and its consequences in combinatorial
  optimization.
\newblock {\em Combinatorica}, 1:169--197, 1981.

\bibitem{Heinonen01}
J.~Heinonen.
\newblock {\em Lectures on analysis on metric spaces}.
\newblock Universitext. Springer-Verlag, New York, 2001.

\bibitem{jl}
W.~B. Johnson and J.~Lindenstrauss.
\newblock Extensions of {L}ipschitz mappings into a {H}ilbert space.
\newblock In {\em Conference in modern analysis and probability (New Haven,
  Conn., 1982)}, pages 189--206. Amer. Math. Soc., Providence, RI, 1984.

\bibitem{KV04}
S.~Khot and N.~Vishnoi.
\newblock The unique games conjecture, integrality gap for cut problems and
  embeddability of negative type metrics into $\ell_1$.
\newblock In {\em Proceedings of the 46th Annual IEEE Conference on Foundations
  of Computer Science}, 2005.
\newblock To appear.

\bibitem{KLMN04}
R.~Krauthgamer, J.~R. Lee, M.~Mendel, and A.~Naor.
\newblock Measured descent: A new embedding method for finite metrics.
\newblock {\em Geom. Funct. Anal.}, 2004.
\newblock To appear.

\bibitem{L04}
J.~R. Lee.
\newblock {O}n distance scales, embeddings, and efficient relaxations of the
  cut cone.
\newblock In {\em Proceedings of the 16th Annual ACM-SIAM Symposium on Discrete
  Algorithms}, pages 92--101, Vancouver, 2005.

\bibitem{LMN04}
J.~R. Lee, M.~Mendel, and A.~Naor.
\newblock Metric structures in {$L_1$}: Dimension, snowflakes, and average
  distortion.
\newblock {\em European J. Combin.}, 2004.
\newblock To appear.

\bibitem{LNdiamond}
J.~R. Lee and A.~Naor.
\newblock Embedding the diamond graph in ${L}_p$ and dimension reduction in
  ${L}_1$.
\newblock {\em Geom. Funct. Anal.}, 14(4):745--747, 2004.

\bibitem{LN04}
J.~R. Lee and A.~Naor.
\newblock Extending {L}ipschitz functions via random metric partitions.
\newblock {\em Invent. Math.}, 160(1):59--95, 2005.

\bibitem{LR99}
T.~Leighton and S.~Rao.
\newblock Multicommodity max-flow min-cut theorems and their use in designing
  approximation algorithms.
\newblock {\em J. ACM}, 46(6):787--832, 1999.

\bibitem{LLR95}
N.~Linial, E.~London, and Y.~Rabinovich.
\newblock The geometry of graphs and some of its algorithmic applications.
\newblock {\em Combinatorica}, 15(2):215--245, 1995.

\bibitem{LinialSaks}
N.~Linial and M.~Saks.
\newblock Low diameter graph decompositions.
\newblock {\em Combinatorica}, 13(4):441--454, 1993.

\bibitem{Mat01}
J.~Matou{\v{s}}ek.
\newblock {\em Lectures on discrete geometry}, volume 212 of {\em Graduate
  Texts in Mathematics}.
\newblock Springer-Verlag, New York, 2002.

\bibitem{SM90}
D.~W. Matula and F.~Shahrokhi.
\newblock The maximum concurrent flow problem.
\newblock {\em J. ACM}, 37(2):318--334, 1990.

\bibitem{MN04}
M.~Mendel and A.~Naor.
\newblock Euclidean quotients of finite metric spaces.
\newblock {\em Adv. Math.}, 189(2):451--494, 2004.

\bibitem{MS86}
V.~D. Milman and G.~Schechtman.
\newblock {\em Asymptotic theory of finite-dimensional normed spaces}, volume
  1200 of {\em Lecture Notes in Mathematics}.
\newblock Springer-Verlag, Berlin, 1986.
\newblock With an appendix by M. Gromov.

\bibitem{NPSS04}
A.~Naor, Y.~Peres, O.~Schramm, and S.~Sheffield.
\newblock Markov chains in smooth normed spaces and {G}romov hyperbolic metric
  spaces.
\newblock Preprint, 2004.

\bibitem{NRS04}
A.~Naor, Y.~Rabani, and A.~Sinclair.
\newblock Quasisymmetric embeddings, the observable diameter, and expansion
  properties of graphs.
\newblock {\em J. Funct. Anal.}
\newblock To appear.

\bibitem{NS02}
A.~Naor and G.~Schechtman.
\newblock Remarks on non linear type and {P}isier's inequality.
\newblock {\em J. Reine Angew. Math.}, 552:213--236, 2002.

\bibitem{Rao99}
S.~Rao.
\newblock Small distortion and volume preserving embeddings for planar and
  {E}uclidean metrics.
\newblock In {\em Proceedings of the 15th Annual Symposium on Computational
  Geometry}, pages 300--306, New York, 1999.

\bibitem{Shmoys95}
D.~B. Shmoys.
\newblock Cut problems and their application to divide-and-conquer.
\newblock In {\em Approximation Algorithms for NP-hard Problems, (D.S.
  Hochbaum, ed.)}, pages 192--235. PWS, 1997.

\bibitem{Vazirani01}
V.~V. Vazirani.
\newblock {\em Approximation algorithms}.
\newblock Springer-Verlag, Berlin, 2001.

\end{thebibliography}

\def\cprime{$'$}

\end{document}